\documentclass{article}
\usepackage{amsmath}
\usepackage[margin=20truemm]{geometry}
\title{\huge Simultaneous diophantine equations of\\
$\sum_{k=1}^{n} x_k^4=\sum_{k=1}^{n} y_k^4,\ \prod_{k=1}^n x_k=\prod_{k=1}^n y_k$, \\
$n=3,4,5,6$}
\author{Seiji Tomita and  Oliver Couto}
\date{}
\begin{document}
\maketitle
\begin{abstract}
In this paper, we proved that there are infinitely many  parametric solutions of $\sum_{k=1}^{n} x_k^4=\sum_{k=1}^{n} y_k^4,\ \prod_{k=1}^n x_k=\prod_{k=1}^n y_k$ where $n=3,4,5,6.$\\
\end{abstract}
\newpage

\centerline{\huge1\quad Introduction}
\vskip\baselineskip

Euler\cite{a} proposed the diophantine equation $a^4+b^4 = c^4+d^4$ for the first time and gave several parametric solutions.
Tariste\cite{a} gave a parametric solution $a^4 + b^4 + c^4 = d^4 + e^4 + f^4$ where $a+b=c$ and $d+e=f$ below.
$(17x-63y)^4 + (23x+57y)^4 + (40x-6y)^4 = (17x+63y)^4 + (23x-57y)^4 + (40x+6y)^4$. \\
Stephane Vandemergel's\cite{c} found three solutions $(a,b,c,d,e,f)=(29,66,124,22,93,116),(54,61,196,28,122,189),\\
(19,217,657,9,511,589)$.
He notes that $r^n+s^n=u^n+v^n$, then $(ru,su,v^2,rv,sv,u^2)$ is a solution, which shows that there are infinitely many solutions for $n<=4$. \\
According to Lander, Parkin, and Selfridge\cite{b} conjecture, we can expect non trivial solution of $\sum_{k=1}^{n} x_k^4=\sum_{k=1}^{n} y_k^4$ for $n>=2$.\\
However there are few parametric solutions for $n>=4$.
In 2016, Tomita\cite{d} found the solution  for $n=3,4$.
\vskip\baselineskip
Using condition $\sum_{k=1}^{n} x_k=\sum_{k=1}^{n} y_k$, we can easily obtain the parametric solution, but as far as I know, Vandemergel is the only one who has used condition $\prod_{k=1}^n x_k=\prod_{k=1}^n y_k$.\\
Therefore, we decided to try to extend $n$ to $6$ using condition $\prod_{k=1}^n x_k=\prod_{k=1}^n y_k$.
\vskip\baselineskip
In this paper, we proved that there are infinitely many  parametric solutions of $\sum_{k=1}^{n} x_k^4=\sum_{k=1}^{n} y_k^4$,\\
$\ \prod_{k=1}^n x_k=\prod_{k=1}^n y_k$ where $n=3,4,5,6.$\\
Furthermore, we showed the parametric solutions of $\sum_{k=1}^{n} x_k^4=\sum_{k=1}^{n} y_k^4$ \ , $\sum_{k=1}^{n} x_k=\sum_{k=1}^{n} y_k$ where $n=3,4,5,6$ in the appendix.
\newpage
\centerline{\huge$2\quad \{x_1^4 + x_2^4 + x_3^4  = y_1^4 + y_2^4 + y_3^4 ,x_1 \cdot x_2 \cdot x_3=y_1 \cdot y_2 \cdot y_3\}$}
\vskip\baselineskip

\textbf{Lemma 1}
Diophantine equation $p^4+q^4=r^4+w^4$ has an infinitely many parametric solutions.

Proof.
\begin{equation}
p^4+q^4=r^4+w^4
\end{equation}

Let $p=kt+u, q=nt+v, r=kt-v, w=nt+u$, then equation $(1)$ becomes to below equation.
\begin{equation}
(4vk^3+4vn^3-4un^3+4uk^3)t^2+(-6u^2n^2-6v^2k^2+6v^2n^2+6u^2k^2)t-4u^3n+4u^3k+4v^3k+4v^3n=0
\end{equation}

Since this is a quadratic equation in $t$, for t to be rational number, the discriminant of the equation must be square number.

Let $U=k/n$ then we obtain
\begin{align}
V^2 &= (-28u^4-28v^4-64vu^3-72v^2u^2-64v^3u)U^4 \nonumber \\
    &+(-64v^4+64u^4+64vu^3-64v^3u)U^3 \nonumber \\
    &+(-72u^4+144v^2u^2-72v^4)U^2 \nonumber \\
    &+(-64vu^3+64v^3u-64v^4+64u^4)U \nonumber \\
    &+64vu^3-72v^2u^2+64v^3u-28u^4-28v^4
\end{align}
Quartic equation $(3)$ has a rational point $Q(U,V)=( \frac{u^3-v^3}{u^3+v^3}, \frac{24u^3v^3(u-v)}{(u^2-uv+v^2)^2(u+v)} )$, then this quartic equation is birationally equivalent to an elliptic curve below.
\begin{align*}
&Y^2-8/3(2v^8-2v^7u+7v^6u^2-5v^5u^3-8u^4v^4-5u^5v^3+7u^6v^2-2u^7v+2u^8)YX/(v^2u^2(u^2-uv+v^2)) \\
&-2304v^3u^3(v^6+uv^5-u^2v^4-2u^3v^3-u^4v^2+u^5v+u^6)Y/((u^2-uv+v^2)^3) \\
&= X^3-32/9(2v^{16}-4v^{15}u+16v^{14}u^2-24u^3v^{13}+32u^4v^{12}-2u^5v^{11}-6u^6v^{10}+41u^7v^9-102u^8v^8 \\
&+41u^9v^7-6u^{10}v^6-2u^{11}v^5+32u^{12}v^4-24u^{13}v^3+16u^{14}v^2-4u^{15}v+2u^{16})X^2/(u^4v^4(u^2-uv+v^2)^2) \\
&+9216(7v^4-12v^3u+10v^2u^2-12vu^3+7u^4)u^6v^6X/((u^2-uv+v^2)^4) \\
&-32768u^2v^2(14v^{20}-52uv^{19}+180u^2v^{18}-424v^{17}u^3+734u^4v^{16}-858u^5v^{15}+702v^{14}u^6-213u^7v^{13} \\
&-1018u^8v^{12}+1979v^{11}u^9-2088u^{10}v^{10}+1979u^{11}v^9-1018v^8u^{12}-213u^{13}v^7+702u^{14}v^6-858v^5u^{15} \\
&+734u^{16}v^4-424u^{17}v^3+180v^2u^{18}-52u^{19}v+14u^{20})/((u^2-uv+v^2)^6)
\end{align*}

The point corresponding to point $Q$ is \\
\begin{align*}
&P(X,Y)=( 32/9(2v^{16}-4v^{15}u+16v^{14}u^2-24u^3v^{13}+32u^4v^{12}-2u^5v^{11}-6u^6v^{10}+41u^7v^9 \\
       &-102u^8v^8+41u^9v^7-6u^{10}v^6-2u^{11}v^5+32u^{12}v^4-24u^{13}v^3+16u^{14}v^2-4u^{15}v+2u^{16}) \\
       &/(u^4v^4(u^2-uv+v^2)^2), \\
       &512/27(2v^{18}+15v^{16}u^2+12u^4v^{14}-43u^6v^{12}+6u^8v^{10}+6u^{10}v^8-43u^{12}v^6+12u^{14}v^4+15u^{16}v^2+2u^{18}) \\
       &/(v^6u^6) ).
\end{align*}

This point $P$ is of infinite order, and the multiples $mP, m = 2, 3$, ...give an infinitely many points.

Hence we can obtain an infinitely many parametric solutions for equation $(1)$.
The proof is completed.
\vskip3\baselineskip

\textbf{Theorem 2}
Simultaneous diophantine equation $\{x_1^4 + x_2^4 + x_3^4  = y_1^4 + y_2^4 + y_3^4 ,x_1 \cdot x_2 \cdot x_3=y_1 \cdot y_2 \cdot y_3\}$ has an infinitely many parametric solutions.
\vskip\baselineskip

Proof.
According to Stephane Vandemergel's\cite{c}, simultaneous diophantine equation $\{x_1^4 + x_2^4 + x_3^4  = y_1^4 + y_2^4 + y_3^4 ,x_1 \cdot x_2 \cdot x_3=y_1 \cdot y_2 \cdot y_3\}$  has a follwing parametric solution.

\begin{align*}
&(x_1,x_2,x_3)=(rp, wp, q^2) \\
&(y_1,y_2,y_3)=(rq, wq, p^2) \\
&Condition: p^4 + q^ 4 = r^4 + w^4.
\end{align*}

From \textbf{Lemma 1}, diophantine equation $p^4+q^4=r^4+w^4$ has an infinitely many parametric solutions, then
simultaneous equation $\{x_1^4+x_2^4+x_3^4 = y_1^4+y_2^4+y_3^4 ,x_1 \cdot x_2 \cdot x_3=y_1 \cdot y_2 \cdot y_3\}$ has an infinitely many parametric solutions.

The proof is completed.
\vskip\baselineskip

\textbf{Example}

From \textbf{Lemma 1}, take $(k,n,t)=(u^3-v^3,u^3+v^3,\frac{3u^2v^2}{u^6-2u^4v^2-2v^4u^2+v^6})$ and $v=1$, then we obtain 
\vskip\baselineskip
\begin{align*}
p = u^7+u^5-2u^3-3u^2+u \\
q = 3u^5+u^2+u^6-2u^4+1 \\
r = 3u^5-u^2-u^6+2u^4-1 \\
s = u^7+u^5-2u^3+3u^2+u
\end{align*}

Hence we obtain a parametric solution.

\begin{align*}
&x_1 = (u^6-3u^5-2u^4+u^2+1)(u^6+u^4-2u^2-3u+1)u \\
&x_2 = (u^6+u^4-2u^2+3u+1)(u^6+u^4-2u^2-3u+1)u^2 \\
&x_3 = (u^6+3u^5-2u^4+u^2+1)^2 \\
&y_1 = (u^6-3u^5-2u^4+u^2+1)(u^6+3u^5-2u^4+u^2+1) \\
&y_2 = (u^6+u^4-2u^2+3u+1)(u^6+3u^5-2u^4+u^2+1)u \\
&y_3 = (u^6+u^4-2u^2-3u+1)^2u^2
\end{align*}

\newpage
\centerline{\huge $3\quad \{x_1^4 + x_2^4 + x_3^4 + x_4^4 = y_1^4 + y_2^4 + y_3^4 + y_4^4,$} 
\centerline{\huge $          x_1 \cdot x_2 \cdot x_3 \cdot x_4 = y_1 \cdot y_2 \cdot y_3 \cdot y_4$ \}}
          
\vskip\baselineskip

\textbf{Theorem 3}
There is a  parametric solution of $x_1^4 + x_2^4 + x_3^4 + x_4^4 = y_1^4 + y_2^4 + y_3^4 + y_4^4.$
\begin{align*}
&x_1 = p(m + r^4 - q^4)    \\      
&x_2 = q(m + p^4 - w^4)    \\      
&x_3 =-r(-m + p^4 - w^4)   \\      
&x_4 = w(m - r^4 + q^4)    \\      
&y_1 = p(m - r^4 + q^4)    \\ 
&y_2 =-q(-m + p^4 - w^4)   \\ 
&y_3 = r(m + p^4 - w^4)    \\ 
&y_4 = w(m + r^4 - q^4)    \\ 
&where\ x_1 \cdot x_2 \cdot x_3 \cdot x_4 = y_1 \cdot y_2 \cdot y_3 \cdot y_4\ and\ p^4+q^4=r^4+w^4. \\
&m\ is\ arbitrary\ integer. 
\end{align*}

Proof.
\begin{equation}
x_1^4 + x_2^4 + x_3^4 + x_4^4 = y_1^4 + y_2^4 + y_3^4 + y_4^4 
\end{equation}
\begin{equation}
x_1 \cdot x_2 \cdot x_3 \cdot x_4 = y_1 \cdot y_2 \cdot y_3 \cdot y_4
\end{equation}

Substitute $(x_1,x_2,x_3,x_4,y_1,y_2,y_3,y_4)=(pm+s, qm+t, rm+u, wm+v, pm-s, qm-t, rm-u, wm-v)$\\
to $(5)$, we obtain
\begin{equation}
2m(m^2pquw+m^2pqrv+m^2ptrw+m^2sqrw+ptuv+squv+stuw+strv)=0
\end{equation}

Take $t= -qu/r, s= -pv/w$ to eliminate equation $(6)$.

Substitute $s$ and $t$ to $(4)$, we obtain
\begin{align}
&-8(q^4w^3r^2u-w^6r^3v+p^4r^3w^2v-w^3r^6u)m^3   \\
&-8(q^4w^3u^3-w^4r^3v^3+p^4r^3v^3-w^3r^4u^3)m=0 \nonumber
\end{align}

Equating to zero the coefficient of $m^3$ in $(7)$, then we obtain
\vskip\baselineskip
$$u = -r(-w^4+p^4), v = w(q^4-r^4)$$  
\vskip\baselineskip
Finally, $(4)$ becomes to as follows,  
$$8m(q-r)(q+r)(q^2+r^2)(p-w)(p+w)(p^2+w^2)(q^4+p^4-w^4-r^4)(-q^4+p^4-w^4+r^4)=0$$
Hence, if $p^4+q^4=r^4+w^4$ then $(4)$ becomes to zero and obtain a parametric solution.
The proof is completed.
\newpage
\centerline{\huge $4\quad \{x_1^4 + x_2^4 + x_3^4 + x_4^4 = y_1^4 + y_2^4 + y_3^4 + y_4^4,$} 
\centerline{\huge $          x_1 \cdot x_2 \cdot x_3 \cdot x_4 = y_1 \cdot y_2 \cdot y_3 \cdot y_4$ \}}
\vskip\baselineskip

\textbf{Theorem 4}
Simultaneous diophantine equation $\{x_1^4 + x_2^4 + x_3^4 + x_4^4 = y_1^4 + y_2^4 + y_3^4 + y_4^4,x_1 \cdot x_2 \cdot x_3 \cdot x_4 = y_1 \cdot y_2 \cdot y_3 \cdot y_4\}$ has an infinitely many parametric solutions.
\vskip\baselineskip

Proof.

According to \textbf{Theorem 3}, simultaneous diophantine equation $\{x_1^4 + x_2^4 + x_3^4 + x_4^4 = y_1^4 + y_2^4 + y_3^4 + y_4^4,x_1 \cdot x_2 \cdot x_3 \cdot x_4 = y_1 \cdot y_2 \cdot y_3 \cdot y_4\}$  has a follwing parametric solution.

\begin{align*}
&x_1 = p(m + r^4 - q^4)    \\      
&x_2 = q(m + p^4 - w^4)    \\      
&x_3 =-r(-m + p^4 - w^4)   \\      
&x_4 = w(m - r^4 + q^4)    \\      
&y_1 = p(m - r^4 + q^4)    \\ 
&y_2 =-q(-m + p^4 - w^4)   \\ 
&y_3 = r(m + p^4 - w^4)    \\ 
&y_4 = w(m + r^4 - q^4)    \\ 
&where\ x_1 \cdot x_2 \cdot x_3 \cdot x_4 = y_1 \cdot y_2 \cdot y_3 \cdot y_4\ and\ p^4+q^4=r^4+w^4. \\
&m\ is\ arbitrary\ integer. 
\end{align*}

From \textbf{Lemma 1}, diophantine equation $p^4+q^4=r^4+w^4$ has an infinitely many parametric solutions, then
simultaneous equation $\{x_1^4 + x_2^4 + x_3^4 + x_4^4 = y_1^4 + y_2^4 + y_3^4 + y_4^4,x_1 \cdot x_2 \cdot x_3 \cdot x_4 = y_1 \cdot y_2 \cdot y_3 \cdot y_4\}$ has an an infinitely many parametric solutions. \\
The proof is completed.
\vskip\baselineskip

\textbf{Example}

Case of $Q(U,V)=( (u^3-v^3)/(u^3+v^3), 24u^3v^3(u-v)/((u^2-uv+v^2)^2(u+v)) )$.
To simplify the result, let $v=1\ and\ m=1$.

\begin{align*}
&p = u^7+u^5-2u^3-3u^2+u \\
&q = u^6+3u^5+u^2-2u^4+1 \\
&r = -u^6+3u^5-u^2+2u^4-1 \\
&w = u^7+u^5-2u^3+3u^2+u
\end{align*}
\vskip\baselineskip

\begin{align*}
&x_1 = u(u^4-3u+u^6-2u^2+1)f1 \\
&x_2 = (3u^5+u^2+u^6-2u^4+1)f1 \\
&x_3 = (-3u^5+u^2+u^6-2u^4+1)f2 \\
&x_4 = u(u^4+3u+u^6-2u^2+1)f2 \\
&y_1 = u(u^4-3u+u^6-2u^2+1)f2 \\
&y_2 = (3u^5+u^2+u^6-2u^4+1)f2 \\
&y_3 = (-3u^5+u^2+u^6-2u^4+1)f1 \\
&y_4 = u(u^4+3u+u^6-2u^2+1)f1 \\
&f1=-1+24u^5-72u^9-192u^{11}+72u^7+24u^{23}-72u^{19}+72u^{21}-192u^{17}+288u^{15}+288u^{13} \\
&f2=1+24u^5-72u^9-192u^{11}+72u^7+24u^{23}-72u^{19}+72u^{21}-192u^{17}+288u^{15}+288u^{13}  \\
&u\ is\ arbitrary.
\end{align*}

\newpage
\centerline{\huge $5\quad \{x_1^4 + x_2^4 + x_3^4 + x_4^4 + x_5^4 = y_1^4 + y_2^4 + y_3^4 + y_4^4 + y_5^4,$} 
\centerline{\huge $         x_1 \cdot x_2 \cdot x_3 \cdot x_4 \cdot x_5=y_1 \cdot y_2 \cdot y_3 \cdot y_4 \cdot y_5$ \}}

\vskip\baselineskip

\textbf{Theorem 5}
Simultaneous diophantine equation $\{x_1^4 + x_2^4 + x_3^4 + x_4^4 + x_5^4 = y_1^4 + y_2^4 + y_3^4 + y_4^4 + y_5^4,\ x_1 \cdot x_2 \cdot x_3 \cdot x_4 \cdot x_5=y_1 \cdot y_2 \cdot y_3 \cdot y_4 \cdot y_5\}$ has an an infinitely many parametric solutions.
\vskip\baselineskip

Proof.
\begin{equation}
x_1^4 + x_2^4 + x_3^4 + x_4^4 + x_5^4 = y_1^4 + y_2^4 + y_3^4 + y_4^4 + y_5^4
\end{equation}
\begin{equation}
 x_1 \cdot x_2 \cdot x_3 \cdot x_4 \cdot x_5=y_1 \cdot y_2 \cdot y_3 \cdot y_4 \cdot y_5
\end{equation}

Substitute $(x_1,x_2,x_3,x_4,x_5)=(kp, mp, q^2, np, rq)$ and $(y_1,y_2,x_3,y_4,y_5)=(kq, mq, p^2, nq, rp)$ to $(8)$, we obtain
$$-(p-q)(p+q)(p^2+q^2)(p^4-n^4-m^4-k^4+r^4+q^4)$$

One of parametric solution of $p^4+q^4+r^4 = k^4+m^4+n^4$ is given by  \textbf{Thorem 2} as follow.

\begin{align*}
&p = (u^6+3u^5-2u^4+u^2+1)^2 \\
&q = (u^6+u^4-2u^2+3u+1)(u^6+u^4-2u^2-3u+1)u^2 \\
&r = (u^6-3u^5-2u^4+u^2+1)(u^6+u^4-2u^2-3u+1)u \\
&k = -(u^6-3u^5-2u^4+u^2+1)(u^6+3u^5-2u^4+u^2+1) \\
&m = (u^6+u^4-2u^2-3u+1)^2u^2 \\
&n = (u^6+u^4-2u^2+3u+1)(u^6+3u^5-2u^4+u^2+1)u \\
\end{align*}
Hence we obtain one of parametric solution of $(8)$.

\begin{align*}
&x_1 = (u^6-3u^5-2u^4+u^2+1)(u^6+3u^5-2u^4+u^2+1)^3 \\
&x_2 = (u^6+u^4-2u^2-3u+1)^2u^2(u^6+3u^5-2u^4+u^2+1)^2 \\
&x_3 = (u^6+u^4-2u^2+3u+1)^2(u^6+u^4-2u^2-3u+1)^2u^4 \\
&x_4 = (u^6+u^4-2u^2+3u+1)(u^6+3u^5-2u^4+u^2+1)^3u \\
&x_5 = (u^6-3u^5-2u^4+u^2+1)(u^6+u^4-2u^2-3u+1)^2u^3(u^6+u^4-2u^2+3u+1) \\
&y_1 = (u^6-3u^5-2u^4+u^2+1)(u^6+3u^5-2u^4+u^2+1)(u^6+u^4-2u^2+3u+1)(u^6+u^4-2u^2-3u+1)u^2 \\
&y_2 = (u^6+u^4-2u^2-3u+1)^3u^4(u^6+u^4-2u^2+3u+1) \\
&y_3 = (u^6+3u^5-2u^4+u^2+1)^4 \\
&y_4 = (u^6+3u^5-2u^4+u^2+1)u^3(u^6+u^4-2u^2-3u+1) \\
&y_5 = (u^6-3u^5-2u^4+u^2+1)(u^6+u^4-2u^2-3u+1)u(u^6+3u^5-2u^4+u^2+1)^2 \\
&u\ is\ arbitrary.
\end{align*}
Since there are infinitely many parametric solutions of $p^4+q^4+r^4 = k^4+m^4+n^4$, then $(8)$ has an infinitely many parametric solutions.

The proof is completed.
\vskip\baselineskip
\newpage
\centerline{\huge $6\quad \{x_1^4 + x_2^4 + x_3^4 + x_4^4 + x_5^4 + x_6^4 = y_1^4 + y_2^4 + y_3^4 + y_4^4 + y_5^4 + y_6^4,$} 
\centerline{\huge $         x_1 \cdot x_2 \cdot x_3 \cdot x_4 \cdot x_5 \cdot x_6=y_1 \cdot y_2 \cdot y_3 \cdot y_4 \cdot y_5 \cdot y_6$ \}}

\vskip\baselineskip

\textbf{Theorem 6}
Simultaneous diophantine equation $\{x_1^4 + x_2^4 + x_3^4 + x_4^4 + x_5^4 + x_6^4 = y_1^4 + y_2^4 + y_3^4 + y_4^4 + y_5^4 + y_6^4,\\ x_1 \cdot x_2 \cdot x_3 \cdot x_4 \cdot x_5 \cdot x_6=y_1 \cdot y_2 \cdot y_3 \cdot y_4 \cdot y_5 \cdot y_6\}$ has an infinitely many parametric solutions.
\vskip\baselineskip

Proof.
\begin{equation}
x_1^4 + x_2^4 + x_3^4 + x_4^4 + x_5^4 + x_6^4 = y_1^4 + y_2^4 + y_3^4 + y_4^4 + y_5^4 + y_6^4
\end{equation}

Substitute $(x_1,x_2,x_3,x_4,x_5,x_6)=(mp, np, q^2, tr, wr, s^2)$ and $(y_1,y_2,y_3,y_4,y_5,y_6)=(mq, nq, p^2, ts, ws, r^2)$ to $(10)$, we obtain
$$x_1^4 + x_2^4 + x_3^4  -( y_1^4 + y_2^4 + y_3^4)=(p-q)(p+q)(p^2+q^2)(-p^4+n^4+m^4-q^4)$$
$$x_4^4 + x_5^4 + x_6^4  -( y_4^4 + y_5^4 + y_6^4)=(s-r)(s+r)(s^2+r^2)(s^4-t^4-w^4+r^4)$$

One of parametric solution of $p^4+q^4 = m^4+n^4$ is given by \textbf{Lemma 1} as follow.

\begin{align*}
&p = u^7+u^5-2u^3-3u^2+u \\
&q = 3u^5+u^2+u^6-2u^4+1 \\
&m = 3u^5-u^2-u^6+2u^4-1 \\
&n = u^7+u^5-2u^3+3u^2+u \\
\end{align*}

Similarly, we obtain
\begin{align*}
&r = v^7+v^5-2v^3-3v^2+v \\
&s = 3v^5+v^2+v^6-2v^4+1 \\
&t = 3v^5-v^2-v^6+2v^4-1 \\
&w = v^7+v^5-2v^3+3v^2+v \\
\end{align*}

Hence we obtain one of parametric solution of $(10)$.

\begin{align*}
&x_1 = (3u^5-u^2-u^6+2u^4-1)(u^7+u^5-2u^3-3u^2+u) \\
&x_2 = (u^7+u^5-2u^3+3u^2+u)(u^7+u^5-2u^3-3u^2+u) \\
&x_3 = (3u^5+u^2+u^6-2u^4+1)^2 \\
&x_4 = (3v^5-v^2-v^6+2v^4-1)(v^7+v^5-2v^3-3v^2+v) \\
&x_5 = (v^7+v^5-2v^3+3v^2+v)(v^7+v^5-2v^3-3v^2+v) \\
&x_6 = (3v^5+v^2+v^6-2v^4+1)^2 \\
&y_1 = (3u^5-u^2-u^6+2u^4-1)(3u^5+u^2+u^6-2u^4+1) \\
&y_2 = (u^7+u^5-2u^3+3u^2+u)(3u^5+u^2+u^6-2u^4+1) \\
&y_3 = (u^7+u^5-2u^3-3u^2+u)^2 \\
&y_4 = (3v^5-v^2-v^6+2v^4-1)(3v^5+v^2+v^6-2v^4+1) \\
&y_5 = (v^7+v^5-2v^3+3v^2+v)(3v^5+v^2+v^6-2v^4+1) \\
&y_6 = (v^7+v^5-2v^3-3v^2+v)^2 \\
&u\ is\ arbitrary.
\end{align*}
Since there are infinitely many parametric solutions of $p^4+q^4 = m^4+n^4$, then $(10)$ has an infinitely many parametric solutions.
The proof is completed.
\vskip\baselineskip

Finally we state conjecture. \\
As mentioned in the appendix, we can find the parametric solution of  $\sum_{k=1}^{n} x_k^4=\sum_{k=1}^{n} y_k^4$ where the condition $\sum_{k=1}^{n} x_k=\sum_{k=1}^{n} y_k$ for any $n$ where $n>=3.$ \\
The same result can be expected for $\prod_{k=1}^n x_k=\prod_{k=1}^n y_k$  as for $\sum_{k=1}^{n} x_k=\sum_{k=1}^{n} y_k$.
\vskip\baselineskip
\textbf{Conjecture}\\  
Simultaneous diophantine equation $\sum_{k=1}^{n} x_k^4=\sum_{k=1}^{n} y_k^4$,\ $\prod_{k=1}^n x_k=\prod_{k=1}^n y_k$ has infinitely many parametric solutions for any $n$ where $n>=3.$
\newpage

\centerline{\huge $ Appendix$} 
\vskip\baselineskip

We show only four parametric solutions $\sum_{k=1}^{n} x_k^4=\sum_{k=1}^{n} y_k^4$ where the condition $x_1 + x_2 + \cdots x_n=y_1 + y_2 + \cdots y_n.$
Using this condition, we can find the parametric solutions for any $n$ where $n>=3.$
\vskip\baselineskip

\textbf{Theorem}  
Simultaneous diophantine equation $\sum_{k=1}^{n} x_k^4=\sum_{k=1}^{n} y_k^4$,\ $\sum_{k=1}^{n} x_k=\sum_{k=1}^{n} y_k$ has an infinitely many parametric solutions for any $n$ where $n>=3.$
\vskip\baselineskip

We prove it for $n=3$.\\
Let $x_1=pt+a, x_2=qt+b, x_3=t, y_1=pt, y_2=qt+a, y_3=t+b$, then we obtain
$$(4ap^3+4bq^3-4aq^3-4b)t^2+(6b^2q^2-6b^2+6a^2p^2-6a^2q^2)t+(4b^3q+4a^3p-4b^3-4a^3q) = 0$$
Since this is a quadratic equation in $t$, for t to be rational number, the discriminant of the equation must be square number.\\
\begin{align*}
V^2 &= -28a^4p^4+(64qa^4+64ab^3-64ab^3q)p^3 \\
&+(-72a^2b^2+72a^2b^2q^2-72q^2a^4)p^2 \\
&+(64a^3b-64a^3bq^3+64q^3a^4)p \\
&+72a^2b^2q^2-72a^2q^4b^2+64b^3q^4a-64b^3aq^3 \\
&-72b^4q^2+64b^4q-28b^4-28b^4q^4+64b^4q^3-28a^4q^4-64a^3qb+64a^3q^4b
\end{align*}

Since this quartic equation has a point \\
$Q(p,V)=((-b^3q+b^3+a^3q)/(a^3), 6b^2(q-1)(a-b)(a^3+a^3q+a^2b-a^2qb+ab^2-ab^2q-b^3q+b^3)/(a^4))$,\\
then this quartic equation is birationally equivalent to an elliptic curve E.\\
Since $E$ has very large terms, we omit it.\\
Based on point $P(p,V)$ corresponding to $Q(p,V)$, we can obtain an infinitely many parametric solutions using group law.

\vskip2\baselineskip

\leftline{\huge $x_1^4 + x_2^4 + x_3^4  = y_1^4 + y_2^4 + y_3^4$}

\begin{flushleft}
$x_1 = 2p^8+6p^7+14p^6-18p^5-19p^4-21p^3+359p^2+111p-42$ \\
$x_2 = 5p^7+15p^6+35p^5-45p^4-58p^3-84p^2+236p-300 $ \\
$x_3 = 2p^7+6p^6+14p^5-18p^4-82p^3-210p^2-82p+174 $ \\
$y_1 = 2(p^7+3p^6+7p^5-9p^4-41p^3-105p^2-41p+87)p $ \\
$y_2 = 4p^7+12p^6+28p^5+27p^4+25p^3+21p^2-227p+306 $ \\
$y_3 = 7p^7+21p^6+49p^5-63p^4-140p^3-294p^2+154p-126 $ \\
$p$\ is\ arbitrary\ integer. 
\end{flushleft}

\leftline{\huge $x_1^4 + x_2^4 + x_3^4 + x_4^4 = y_1^4 + y_2^4 + y_3^4 + y_4^4$ }

\begin{flushleft}
$x_1 = 2p^9+2p^8+2p^7+38p^6+101p^5+101p^4+317p^3+1073p^2+1416p+1848$  \\      
$x_2 = 5p^8+5p^7+5p^6+95p^5+242p^4+242p^3+194p^2+1076p+2448$   \\      
$x_3 = 2p^8+2p^7+2p^6-25p^5-25p^4-25p^3-565p^2-1594p-2280$ \\      
$x_4 = 2p^8+2p^7+2p^6+38p^5+38p^4+38p^3+254p^2+254p-432$  \\      
$y_1 = 2*(p^8+p^7+p^6+19p^5+19p^4+19p^3+127p^2+127p-216)p$  \\ 
$y_2 = 4p^8+4p^7+4p^6+13p^5+13p^4+13p^3-311p^2-1340p-2712$ \\ 
$y_3 = 7p^8+7p^7+7p^6+133p^5+280p^4+280p^3+448p^2+1330p+2016$    \\ 
$y_4 = 2p^8+2p^7+2p^6+101p^5+101p^4+101p^3+1073p^2+2102p+1416$    \\ 
$p$\ is\ arbitrary\ integer. 
\end{flushleft}
\newpage

\leftline{\huge $x_1^4 + x_2^4 + x_3^4 + x_4^4 + x_5^4 = y_1^4 + y_2^4 + y_3^4 + y_4^4 + y_5^4$} 

\begin{flushleft}
$x_1 = 6p^9+9p^8-153p^7-1305p^6-5283p^5-13950p^4-24624p^3-29448p^2-21564p-8232 $ \\
$x_2 = 12p^9+135p^8+621p^7+1017p^6-1809p^5-13500p^4-33066p^3-46170p^2-36828p-14160 $ \\
$x_3 = 6p^9+54p^8+324p^7+1872p^6+6471p^5+14787p^4+21393p^3+20097p^2+10512p+2892 $ \\
$x_4 = 78p^8+564p^7+996p^6-1767p^5-13323p^4-31413p^3-42105p^2-31926p-11784 $ \\
$x_5 = 6p^8+114p^7+546p^6+1230p^5+1014p^4-1578p^3-5286p^2-6150p-2964 $ \\
$y_1 = 6p^9+114p^8+546p^7+1230p^6+1014p^5-1578p^4-5286p^3-6150p^2-2964p $ \\
$y_2 = 84p^8+678p^7+1542p^6-537p^5-12309p^4-32991p^3-47391p^2-38076p-14748 $ \\
$y_3 = 6p^9-3p^8-381p^7-2397p^6-7743p^5-15978p^4-21468p^3-18876p^2-9264p-2304 $ \\
$y_4 = 66p^8+336p^7-96p^6-4227p^5-15351p^4-28257p^3-31533p^2-19626p-5856 $ \\
$y_5 = 12p^9+129p^8+507p^7+471p^6-3039p^5-14514p^4-31488p^3-40884p^2-30678p-11196 $ \\
$p$\ is\ arbitrary\ integer. 
\end{flushleft}
\newpage

\leftline{\huge $x_1^4 + x_2^4 + x_3^4 + x_4^4 + x_5^4 + x_6^4 = y_1^4 + y_2^4 + y_3^4 + y_4^4 + y_5^4 + y_6^4$} 

\begin{flushleft}
$x_1 = 120p^{11}-1838p^{10}+4347p^9-37848p^8+36468p^7-32262p^6+226989p^5+649908p^4+105264p^3 $ \\
       \qquad $+551372p^2-357312p-17472 $ \\
$x_2 = 24p^{11}-96p^{10}-1235p^9-5778p^8-4644p^7+2979p^6-97665p^5+43272p^4-967764p^3+182295p^2 $ \\
       \qquad $-789772p+112992 $ \\
$x_3 = 24p^{11}-562p^{10}+6727p^9-3828p^8+122256p^7-129822p^6-119295p^5-932184p^4-1538640p^3 $ \\
       \qquad $-698048p^2-1246060p+85440 $ \\
$x_4 = 72p^{11}-782p^{10}-1739p^9-22284p^8-31032p^7+78147p^6+309p^5+736344p^4-838884p^3 $ \\
       \qquad $+678125p^2-841228p+112992 $ \\
$x_5 = 24p^{11}-48p^{10}-1921p^9-6282p^8-21150p^7-23409p^6-22497p^5+141246p^4-274692p^3 $ \\
       \qquad $+311175p^2-293942p+61536 $ \\
$x_6 = 48p^{10}-686p^9-504p^8-16506p^7-26388p^6+75168p^5+97974p^4+693072p^3+128880p^2+495830p $ \\
       \qquad $-51456 $ \\
$y_1 = 24p^{11}-686p^{10}+2103p^9-9720p^8+11268p^7+124965p^6+45303p^5+453852p^4-289500p^3 $ \\
       \qquad $+55775p^2-253344p-10080 $ \\
$y_2 = 96p^{10}-1372p^9-1008p^8-33012p^7-52776p^6+150336p^5+195948p^4+1386144p^3+257760p^2 $ \\
       \qquad $+991660p-102912 $ \\
$y_3 = 120p^{11}-1742p^{10}+2975p^9-38856p^8+3456p^7-85038p^6+377325p^5+845856p^4+1491408p^3 $ \\
       \qquad $+809132p^2+634348p-120384 $ \\
$y_4 = 24p^{11}-782p^{10}+3475p^9-8712p^8+44280p^7+177741p^6-105033p^5+257904p^4-1675644p^3 $ \\
       \qquad $-201985p^2-1245004p+92832 $ \\
$y_5 = 72p^{11}-1200p^{10}+5537p^9-20838p^8+79362p^7-81042p^6+53847p^5-141138p^4-716688p^3 $ \\
       \qquad $-73338p^2-801686p+33984 $ \\
$y_6 = 24p^{11}+48p^{10}-3293p^9-7290p^8-54162p^7-76185p^6+127839p^5+337194p^4+1111452p^3 $ \\
       \qquad $+568935p^2+697718p-41376 $ \\
$p$\ is\ arbitrary\ integer. 
\end{flushleft}

\end{document}